\documentclass
{article}
\usepackage[T2A]{fontenc}
\usepackage[cp1251]{inputenc}
\usepackage[english,russian]{babel}
\usepackage{amsmath}
\usepackage{amsfonts,amssymb,mathrsfs,amscd}
\usepackage{verbatim}
\usepackage{eucal}
\usepackage{graphicx}
\usepackage{enumerate}
\usepackage{esint}
\usepackage{upgreek}
\let\No=\textnumero

\usepackage[hyper]{izv2e}
\JournalName{}

\numberwithin{equation}{section}

\theoremstyle{plain}
\newtheorem{theorem}{Теорема}
\newtheorem{maintheorem}{Основная теорема}

\newtheorem{theo}{Теорема} 

\theoremstyle{definition}
\newtheorem{proof}{Доказательство}
\newtheorem{remark}{Замечание}

\renewcommand{\leq}{\leqslant} 
\renewcommand{\geq}{\geqslant}
\newcommand{\RR}{\mathbb{R}} 
\newcommand{\CC}{\mathbb{C}} 
\newcommand{\NN}{\mathbb{N}} 

\newcommand{\DD}{\mathbb{D}} 

\DeclareMathOperator{\rad}{rad}

\DeclareMathOperator{\ord}{ord}

\DeclareMathOperator{\dd}{d}

\begin{document} 

\title{Ограничение снизу субгармонической функции логарифмом модуля целой
}
\author[B.\,N.~Khabibullin]{Б.\,Н.~Хабибуллин}
\address{Башкирский государственный университет,\\
Институт математики с ВЦ УФИЦ РАН}
\email{khabib-bulat@mail.ru} 

\date{22.03.2022}
\udk{517.538 : 517.574}

\maketitle
	
	\begin{fulltext}

\begin{abstract} 
Пусть $u\not\equiv -\infty$ --- субгармоническая функция на комплексной плоскости $\mathbb \CC$. 
Тогда  для любой функции 
$r\colon \mathbb C\to (0,1]$,  удовлетворяющей  условию 
$$
\inf_{z\in \mathbb C}\frac{\ln r(z)}{\ln (2+|z|)}>-\infty,
$$
существует целая функция $f\not\equiv 0$, для которой 
$$
\ln |f(z)|\leq \frac{1}{2\pi}\int_0^{2\pi}u\bigl(z+r(z)e^{i\theta}\bigr)\,{\mathrm d}\theta
\quad\text{при всех $z\in \CC$}.
$$
Аналогичный результат установлен для субгармонических функций конечного порядка с неравенствами вида $\ln |f(z)|\leq u(z)$ во всех точках $z\in \mathbb C\setminus E$, где  исключительное множество $E$  мал\'о в терминах $d$-мерного обхвата Хаусдорфа переменного радиуса $r$.   	


\end{abstract} 

\markright{Ограничение снизу субгармонической функции}

\footnotetext[0]{Исследование выполнено при финансовой поддержке Российского научного фонда, проект №~22-21-00026, 
\href{https://rscf.ru/project/22-21-00026/}{https://rscf.ru/project/22-21-00026/}}


\section{Введение}\label{s10}

Для произвольной субгармонической на $\CC$  функции $u\not\equiv -\infty$, конечно же, не всегда  существует целая функция $f\not\equiv 0$, для которой $\ln|f|\leq u$ на $\CC$ или на подмножестве в $\CC$. Одна из причин этого в том,  что  {\it $(-\infty)$-множество\/} 
\cite[3.5]{Rans}
\begin{equation}\label{fytu}
(-\infty)_u:=\bigl\{z\in \CC\bigm| u(z)=-\infty\bigr\}
\end{equation} 
{\it субгармонической функции\/} $u$ может быть не локально конечным или  даже оказаться всюду плотным в $\CC$. 
При необходимости какого-либо ограничения снизу функции $u$ через $\ln|f|$ с целой функцией $f\not\equiv 0$ возможны различные подходы, из  которых здесь рассматриваются два: неравенства вида $\ln |f(z)|\leq u(z)$  для точек $z$, лежащих  вне некоторого по возможности малого исключительного множества $E\subset \CC$, или  ограничения снизу 
через $\ln |f|$ для интегральных средних  от функции $u$  по кругам  или окружностям переменного малого  радиуса.
Последний вариант  часто даже  предпочтительнее, поскольку, во-первых, он аддитивен в отличие от наиболее часто используемых ограничений снизу для точной верхней грани функции $u$ по кругам, а во-вторых, из него следуют и ограничения снизу вне исключительных множеств $E$, что в форме критерия отражено, например, в нашей статье  \cite[теорема 2]{Kha21}.

Для формулировки основного результата  потребуется  

\section{Некоторые определения, обозначения и соглашения}
 Одноточечные множества $\{x\}$ часто записываем без фигурных скобок, т.е. просто как $x$. Так,  
 $\overline \RR:=-\infty \cup \RR\cup +\infty$ --- {\it расширение\/}
 {\it вещественной оси\/} $\RR$ путём порядкового пополнения $\RR$ двумя точками  
$$
-\infty:=\inf \RR\notin \RR,  \quad +\infty:=\sup \RR\notin \RR,
$$ неравенствами $-\infty\leq x\leq +\infty$ для любого $x\in \overline \RR$  и естественной порядковой топологией. Через 
 $\RR^+:=\bigl\{x\in \RR\bigm| x\geq 0\bigr\}$ обозначаем {\it положительную полуось.\/} По определению 
$\sup \varnothing:=-\infty$ и $\inf \varnothing:=+\infty$ для {\it пустого множества\/} $\varnothing$.
Символом  $0$, кроме нуля,   могут обозначаться {\it нулевые\/} функции,  меры  и пр.

Для $x\in X\subset \overline \RR$ его {\it положительную часть\/} обозначаем как $x^+:=\sup\{0,x \}$. {\it Расширенной числовой функции\/} $f\colon S\to \overline \RR$ сопоставляем её {\it положительную часть\/} $f^+\colon s\underset{s \in S}{\longmapsto} (f(s))^+\in \overline{\RR}^+$. Как обычно, пишем $f\not\equiv c$, если функция $f$ принимает хотя бы одно значение, отличное от $c$. 

Для $x_0\in \RR$ и расширенной числовой функции $m\colon x_0+\RR^+ \to \overline \RR$ определим  
\begin{equation}\label{senu0:a} 
\ord[m]:=\limsup_{x\to +\infty} \frac{\ln \bigl(1+m^+(x)\bigr)}{\ln x}\in
\overline \RR^+:=\RR^+\cup +\infty
 \end{equation} 
{\it --- порядок\/}  (роста) функции $m$ (около $+\infty$) \cite{Boas}, \cite{Levin56}, \cite{Levin96}, 
\cite{Kiselman}, \cite[2.1]{KhaShm19}, а для произвольной 
 функции $u\colon \CC\to \overline \RR$ с  {\it радиальный функцией роста\/}
\begin{equation} 
{\mathrm M}_u \colon r\underset{r\in \RR^+}{\longmapsto}
\sup\bigl\{u(z)\bigm| |z|=r\bigr\}
\label{u}
\end{equation} 
по определению 
$$
\ord[u]\overset{\eqref{senu0:a}}{:=}\ord[{\mathrm M}_u]
$$
--- {\it порядок\/} (роста) функции $u$ \cite{Boas}, \cite{Levin96}, \cite{Kiselman}, \cite[Замечание 2.1]{KhaShm19}.  

{\it Распределением масс\/} называем  {\it положительную   меру Радона\/}  \cite{EG},  \cite[Дополнение A]{Rans}, \cite[гл.~3]{HK}. 
Для распределений масс    {\it на\/} $\CC$, как правило, не указываем, где  они заданы.   
Для  субгармонической  в области из $\CC$ функции $u\not\equiv -\infty$ действие  на неё   {\it оператора Лапласа\/} ${\bigtriangleup}$  в смысле теории обобщённых функций  определяет её {\it распределение масс Рисса\/}
\begin{equation}\label{Riesz}
\varDelta_u:=\frac{1}{2\pi}{\bigtriangleup}u
\end{equation}
 в этой области  \cite{HK}, \cite{Rans}, \cite{Landkof}, \cite{Az}. Через 
\begin{equation*}
\begin{split}
D_z(r)&:=\bigl\{w \in \CC \bigm| |w-z|<r\bigr\}\\
\overline{D}_z(r)&:=\bigl\{w \in \CC \bigm| |w-z|\leq r\bigr\},\\
\partial \overline{D}_z(r)&:=\overline{D}_z(r)\setminus {D}_z(r)
\end{split}
\end{equation*}
обозначаем  соответственно {\it открытый} и   {\it замкнутый круги,\/} а также  {\it окружность  радиуса\/ $r\in \RR^+$
 с центром\/} $z\in \CC$, а   $\DD:=D_0(1)$ и  $\overline \DD:=\overline D_0(1)$, а также  
$\partial \overline \DD:=\partial \overline{D}_0(1)$ --- соответственно {\it открытый\/} и {\it замкнутый 
единичные круги,\/} а также {\it единичная окружность} в  $\CC$.

При $r\in \RR^+$ для распределения  масс   $\mu$  через 
\begin{equation}\label{df:nup} 
\mu_z^{\rad} (r):=\mu \bigl(\,\overline D_z(r)\bigr),\quad \mu^{\rad}(r):=\mu_0^{\rad}(r)
=\mu\bigl(r\overline \DD\bigr)
\end{equation} 
обозначаем  {\it радиальные непрерывные справа считающие  функции  распределения масс  $\mu$ с
центрами\/}  соответственно {\it в точке\/ $z\in \CC$} и {\it в нуле.\/}

{\it Порядок распределения масс\/} $\mu$ определяется через \eqref{df:nup} как 
\begin{equation}\label{ordnu}
\ord[\mu]\overset{\eqref{senu0:a}}{:=}\ord\bigl[\mu^{\rad}\bigr].
\end{equation}
.

Развивая \cite[II]{Carleson}, \cite{Rodgers}, \cite{Federer},    \cite[гл. 2]{EG}, \cite{Eid07}, \cite[5.2]{VolEid13}
\cite[определение 3]{Kha22}, для $d\in  \RR^+$, функции  $r\colon \CC\to   \overline \RR^+\setminus 0$ и {\it гамма-функции\/} $\Gamma$ внешнюю   меру 
 \begin{equation}\label{mr}
{\mathfrak m}_d^r\colon S\underset{S\subset \CC}{\longmapsto}  
\inf \Biggl\{\sum_{k} \dfrac{\pi^{d/2}}{\Gamma (1+d/2)}r_k^d\biggm| S\subset \bigcup_{k} 
\overline D_{z_k}(r_k), \, z_k\in \CC, \, r_k \leq  r(z_k)\Biggr\}
\end{equation}
 называем {\it $d$-мерным обхватом Хаусдорфа переменного радиуса обхвата\/ $r$.\/} При этом  
через {\it постоянные функции\/} $r>0$ определяется {\it $d$-мерная    мера Хаусдорфа\/}
\begin{equation}\label{hH}
{\mathfrak m}_d\colon S\underset{S\subset \CC}{\longmapsto}  \lim_{0<r\to 0} {\mathfrak m}_d^r(S)
\underset{r>0}{\geq} {\mathfrak m}_d^r(S)\geq {\mathfrak m}_d^\infty(S),
\end{equation}
являющаяся регулярной  мерой Бореля. По определению \eqref{mr}, очевидно, 
\begin{equation}\label{mts}
{\mathfrak m}_d\geq {\mathfrak m}_d^{r}\geq {\mathfrak m}_d^{t}\geq {\mathfrak m}_d^{\infty}
\quad\text{для любых пар  функций $r\leq t$.} 
\end{equation}

В частности,   ${\mathfrak m}_2$ ---  это {\it плоская мера  Лебега на\/} $\CC$,  а для любой  липшицевой кривой $L$ в $\CC$  сужение  меры Хаусдорфа  ${\mathfrak m}_1$ на образ этой  кривой в $\CC$ --- это мера длины дуги на этой кривой \cite[3.3.4A]{EG}, в соответствии  с чем используем обозначение  ${\mathfrak m}_1$ и для линейной меры  Лебега на $\RR$.
Кроме того,  $0$-мерная   мера Хаусдорфа ${\mathfrak m}_0$ множества -- это число элементов в этом множестве. При этом для плоскости $\CC$  имеет смысл рассматривать  только $d\in [0,2]$, поскольку 
\begin{equation}\label{d2}
{\mathfrak m}_d={\mathfrak m}_d^r=0\quad \text{при любом    $d>2$.}
\end{equation}

Как и  в \cite[предисловие]{EG}, расширенная числовая функция {\it интегрируема\/} по мере  Радона на множестве, если интеграл от неё по этой мере  корректно определён значением из $\overline \RR$. Интегрируемая функция {\it суммируема,\/}  если соответствующий  интеграл от неё конечен, т.е. принимает значения из $\RR$.

Для функций  $r\colon S\to \RR^+$  на подмножестве $S\subset \CC$ будут использованы  {\it интегральные средние 
 с переменным радиусом $r$  по  окружностям}
\begin{equation}
v^{\circ r}\colon z\underset{z\in S}{\longmapsto}\frac{1}{2\pi r(z)} \int_{\partial D_z(r(z))} 
 v \dd \mathfrak m_1  
\label{vpC}
\end{equation}
для  ${\mathfrak m}_1$-интегрируемых  функций  $v$  на окружностях $\partial D_z\bigl(r(z)\bigr)$ при $z\in S$, а также 
{\it интегральные средние  с переменным радиусом $r$   по кругам\/}
\begin{equation}
v^{\bullet r}\colon z\underset{z\in S}{\longmapsto}
\frac{1}{\pi (r(z))^2} \int_{\overline  D_z(r(z))}  v \dd \mathfrak m_2  
\label{vpD}
\end{equation}
для $\mathfrak m_2$-интегрируемых  функций $v$ на кругах $\overline  D_z\bigl(r(z)\bigr)$ при $z\in S$. 
Для любой функции $u$, субгармонической  на открытой окрестности объединения кругов 
\begin{equation}\label{Sdrcup}
 S^{\cup r}:=\bigcup_{z\in S}\overline  D_z\bigl(r(z)\bigr)\subset \mathbb C 
\end{equation}
из  \eqref{vpD}, имеем неравенства \cite[теорема 2.6.8]{Rans}
\begin{equation}\label{vbc}
u(z)\leq u^{\bullet r}(z)\leq u^{\circ r}(z) \leq \sup_{D_z(r(z))}u
\quad\text{при всех }z\in S^{\cup r}.
\end{equation}

\section{Формулировка основного результата}\label{S8se}

\begin{maintheorem}
Пусть $u\not\equiv -\infty$ --- субгармоническая функция на $\CC$, а функция 
 $r\colon \CC\to (0,1]$ удовлетворяет условию 
\begin{equation}\label{qr}
\inf_{z\in \CC} \frac{\ln r(z)}{\ln(2+ |z|)}>-\infty
\end{equation}
 Тогда существует целая функция $f\not\equiv 0$, для которой 
\begin{align}
\ln |f(z)|\overset{\eqref{vpD}}{\leq} u^{\bullet r}(z) &\overset{\eqref{vbc}}{\leq} u^{\circ r}(z) 
\overset{\eqref{vbc}}{\leq}\sup_{D_z(r(z))}u \quad\text{при всех $z\in \CC$},
\label{ubc}
\\
{\mathrm M}_{\ln |f|}(R)&\overset{\eqref{u}}{\leq}  {\mathrm M}_u\bigl(R+{\mathrm M}_r(R) \bigr)\quad\text{при всех $R\in \RR^+$}.
\label{Mulnf}
\end{align}
Если для  этой субгармонической  функции $u\not\equiv -\infty$ её распределение масс Рисса $\varDelta_u$  конечного порядка $\ord[\varDelta_u]\overset{\eqref{ordnu}}{<}+\infty$,  то для  каждого  $d\overset{\eqref{d2}}{\in} (0,2]$ целую функцию $f\not\equiv 0$, удовлетворяющую\/  \eqref{ubc}--\eqref{Mulnf}, можно подобрать и  так, что 
\begin{align}
\ln |f(z)|&\leq u(z) \quad\text{при всех $z\in \CC\setminus E$, где}
\label{{uE1}u}\\
{\mathfrak m}_d^{r}(E\cap S)&\leq \sup_{z\in S} r(z)
\quad\text{для  любого  $S\subset  \CC$}.
\label{{uE1}E}
\end{align}
\end{maintheorem}

\begin{remark}\label{rem3}
Условие   \eqref{qr} для  $r\colon \CC\to (0,1]$ эквивалентно существованию достаточно малого   $c\in \RR^+\setminus 0$ и больших   
$R\in \RR^+$ и $P\in \RR^+$, для которых 
\begin{equation}\label{nubstrr}
 r(z)\geq\begin{cases}
 c>0&\text{ при всех $z\in  R\overline \DD$},\\ 
\dfrac{1}{(1+|z|)^P}&\text{ при всех $z\in \CC\setminus R\overline \DD$}.
\end{cases}
\end{equation}
В частности, при такой  функции $r$  для исключительного множества $E\subset \CC$ соотношение \eqref{{uE1}E} при выборе $S:=\CC\setminus t\overline \DD$ даёт 
$$
{\mathfrak m}_d^r(E_b\setminus t\overline \DD) \overset{\eqref{mts}}{=}O\Bigl(\frac{1}{t^P}\Bigr)\quad\text{при $t\to +\infty$}.
$$
Отсюда, например,   при выборе $d:=1$ для  любого сколь угодно большого $P>0$  исключительное множество $E$, 
удовлетворяющее \eqref{{uE1}E}, может быть выбрано так, что длина её части, лежащей на (образах) липшицевых кривых  $L\subset \CC$
 не только конечной линейной меры Лебега ${\mathfrak m}_1$,  но и быстро уменьшается по мере ${\mathfrak m}_1$ со степенной  скоростью $t^{-P}$ при $t\to +\infty$ вне  кругов  $t\DD$, т.е. на   $L\setminus t\DD$.  
\end{remark}

\section{Доказавтельство основной теоремы}

 Из  \cite[п.~1.3, следствие  2]{KhaBai16} с комментарием после формулировки, основанном на 
 \cite[п.~2.4, доказательство следствия 2]{KhaBai16},  и в форме из    \cite[п.~4, лемма 5.1]{BaiKhaKha17} и \cite[теорема 1]{Kha21} для любой  
 субгармонической на $\CC$ функции $u\not\equiv -\infty$  и  любого числа $P\in \RR^+$  с соответствующей функцией 
\begin{equation}\label{pP}
p\colon z\underset{z\in \CC}{\longmapsto} \frac{1}{\bigl(1+|z|\bigr)^P}
\end{equation}
найдётся целая функция $f_P\not\equiv 0$,  для которой
$\ln\bigl|f_P(z)\bigr|\leq u^{\bullet p}(z)$
при всех $z\in \CC$. Отсюда для функции $r$ с оценками снизу в форме  \eqref{nubstrr} из замечания  \ref{rem3}  при выборе достаточно большого $R\in \RR^+$ имеем  $r(z)\geq p(z)$ и $\ln\bigl|f_P(z)\bigr|\leq u^{\bullet r}(z)$
при всех $z\in \CC\setminus R\overline \DD$.
В то же время  по неравенству   \eqref{nubstrr} из замечания  \ref{rem3} 
имеем $r(z)\geq c$ при всех  $z\in R\overline \DD$ для некоторого числа $c\in \RR^+\setminus 0$, откуда 
при всех  $z\in \overline \DD$ получаем $u^{\bullet r}(z)\geq u^{\bullet c}(z)$. Правая часть в этом неравенстве непрерывна, следовательно, существует  $c_R\in \RR^+\setminus 0$, для которого $u^{\bullet r}(z)\geq c_R$ при всех  $z\in R\overline \DD$. Отсюда  при достаточно малом значении числа $a\in \RR^+\setminus 0$ для целой функции $f:=af_P$ получаем требуемое  
\eqref{ubc} и, как  следствие принципа  максимума,   \eqref{Mulnf}.

При обосновании заключительной части основной теоремы  используем следующие обозначения.  
Для   $s\in  \mathbb R^+$ и  функции $p\colon \mathbb C \to \mathbb R$ положим 
\begin{equation}\label{Sdrp}
p^{\vee s}(z)\underset{z\in \CC}{:=} \sup \Bigl\{ p(w)\Bigm| w\in \overline D_z(s)\Bigr\}\in \overline{\mathbb R},
\end{equation}
\begin{theorem}\label{lemmuq} Пусть  $\mu$ --- распределение масс на $\CC$, а  $p\colon \CC \to \RR^+\setminus 0$ 
--- борелевская ограниченная функция и 
\begin{equation}\label{r0}
s:=\sup_{z\in \CC} p(z) <+\infty.
\end{equation} 
Тогда при любом  $d\in (0,2]$ множество 
\begin{equation}\label{Er}
 E:= \biggl\{z\in \mathbb C  \biggm| \int_0^{p(z)}\frac{\mu_z^{\rad}(t)}{t}\,{\rm d} t>\frac{1}{d}\biggr\}\subset \mathbb C 
\end{equation}
для каждого $S\subset \CC$ в обозначении \eqref{Sdrp} удовлетворяет  неравенству 
\begin{equation}\label{mmES}
{\mathfrak m}_d^p(E\cap S) \leq 60 \int_{S_s} \bigl(p^d\bigr)^{\vee s} \dd \mu\in \overline \RR^+
\end{equation}
для некоторого борелевского подмножества $S_s\overset{\eqref{Sdrcup}}{\subset} S^{\cup (2s)}$.
\end{theorem}
\begin{proof}
Если для точки $z\in \CC$ выполнено 
\begin{equation}\label{Qp}
\mu_z^{\rad}(t)\leq p^{-d}(z)t^d\quad\text{при всех $t\leq p(z)$},
\end{equation} 
то эта  точка $z$ не принадлежит множеству $E$ из \eqref{Er}, поскольку для неё
\begin{equation*}
\int_0^{p(z)}\frac{\mu_z^{\rad}(t)}{t}\,{\rm d} t\overset{\eqref{Qp}}{\leq} 
p^{-d}(z)\int_0^{p(z)} t^{d-1}\dd t=\frac{1}{d}. 
\end{equation*}
Таким образом, для каждой точки   $z\overset{\eqref{Er}}{\in} E$  отрицание \eqref{Qp} даёт 
\begin{equation}\label{murt}
t_z^d\leq \frac{1}{p^{-d}(z)}\mu_z^{\rad}(t_z)\overset{\eqref{Qp}}{=} p^d(z) \mu_z^{\rad}(t_z)
\quad\text{при  некотором  $t_z\in \bigl(0,p(z)\bigr]$}, 
\end{equation}  
 и система  кругов $\bigl\{\overline  D_z(t_z)\bigr\}_{z\in E}$, очевидно, покрывает 
\begin{equation}\label{DzE}
E\subset \bigcup_{z\in E} \overline D_z(t_z), \quad 0<t_z\leq p(z)\overset{\eqref{r0}}{\leq} s. 
\end{equation}
\begin{theo}[Безиковича о покрытиях {(\cite[лемма  3.2]{Landkof}, \cite[2.8.14]{Federer}, \cite{RKrantz},  \cite{BEr}, \cite{FL}, \cite[I.1, Замечания]{Gusman}, \cite{Sullivan})}] 
Пусть  $t\colon z\underset{z\in E}{\longmapsto} t_z\in {\mathbb R}^+\setminus 0$ --- функция на  $E\subset {\mathbb C}$. Если  $t$ или $E$ ограничены, то  для  некоторого $N\subset \NN$ найдётся  последовательность попарно  различных точек $z_k\in E$, $k\in N\subset {\mathbb N}$,  для которых $E\subset \bigcup\limits_{k\in N} \overline D_{z_k}\bigl(t_{z_k}\bigr)$ и 
пересечение любых  $20$ различных кругов  $\overline D_{z_k}\bigl(t_{z_k}\bigr)$ пусто.
\end{theo}

По теореме Безиковича о покрытиях можно перейти от \eqref{DzE} к не более чем счётной системе кругов 
$\bigl\{\overline  D_{z_k}(t_k)\bigr\}_{z_k\in E}$, покрывающей множество 
\begin{equation}\label{DzEk}
E\subset \bigcup_{z_k\in E} \overline  D_{z_k}(t_k), \quad 0<t_k:=t_{z_k}\leq p(z_k)\overset{\eqref{r0}}{\leq} s. 
\end{equation}

Рассмотрим произвольное   $S\subset \CC$, для которого  множество 
\begin{equation}\label{S}
S_s:=\bigcup\Bigl\{\overline  D_{z_k}( t_k)\Bigm| S\cap \overline D_{z_k}(t_k)\neq \varnothing\Bigr\} 
\overset{\eqref{DzEk}}{\subset} \bigcup_{z\in S}\overline D_z(2s)
\overset{\eqref{Sdrcup}}{=}S^{\cup (2s)},
\end{equation}
очевидно, борелевское,  а также имеем неравенства 
\begin{multline}\label{djkf}
\sum_{S\cap \overline D_{z_k}(t_k)\neq \varnothing} t_{k}^d \overset{\eqref{murt}}{\leq}
\sum_{S\cap \overline D_{z_k}(t_k)\neq \varnothing} 
 p^d(z_k) \mu_{z_k}^{\rad}(t_k)
\\
= \sum_{S\cap \overline D_{z_k}(t_k)\neq \varnothing} \int _{\overline D_{z_k}(t_k)}p^d(z_k)\dd \mu(z).
\end{multline}

{\it Характеристическую функцию множества\/}  $S$ обозначаем через 
\begin{equation}\label{SdrS}
\boldsymbol{1}_S\colon z\underset{z\in \mathbb C}{\longmapsto} \begin{cases}
1&\text{ если $z\in S$},\\
0&\text{ если $z\notin S$}.
\end{cases}
\end{equation}

Ввиду $z_k\overset{\eqref{r0}}{\in} \overline D_z(s)$ для всех $z\in \overline D_{z_k}(t_k)$ 
и соответствующего  неравенства  $p^d(z_k)\overset{\eqref{r0}}{\leq} \bigl(p^d\bigr)^{\vee s}(z)$  для всех таких $z\in \overline D_{z_k}(t_k)$, неравенства \eqref{djkf} дают 
\begin{multline*}
\sum_{S\cap \overline D_{z_k}(t_k)\neq \varnothing} t_{k}^d
\overset{\eqref{Sdrp},\eqref{DzEk}}{\leq} 
\sum_{S\cap \overline D_{z_k}(t_k)\neq \varnothing} \int _{\overline D_{z_k}(t_k)}\bigl(p^d\bigr)^{\vee s}
{\rm \,d} \mu
\\
\overset{\eqref{S},\eqref{SdrS}}{=}
 \sum_{S\cap \overline D_{z_k}(t_k)\neq \varnothing} \int_{S_s}
\boldsymbol{1} _{\overline D_{z_k}(t_k)} \bigl(p^d\bigr)^{\vee s} \dd \mu\\
=  \int_{S_s} \Biggl(\sum_{S\cap \overline D_{z_k}(t_k)\neq \varnothing}\boldsymbol{1} _{\overline D_{z_k}(t_k)}\Biggr) 
\bigl(p^d\bigr)^{\vee s}  \dd \mu.
\end{multline*} 
Отсюда по теореме Безиковича о покрытиях  внутренняя подынтегральная сумма в последнем интеграле не превышает  $19$ и 
\begin{equation}\label{SDyU} 
\sum_{S\cap \overline D_{z_k}(t_k)\neq \varnothing} t_{k}^d\leq 
19\int_{S_s} \bigl(p^d\bigr)^{\vee s} \dd \mu.
\end{equation} 
Кроме того, по определению ${\mathfrak m}_d^p$-обхвата Хаусдорфа \eqref{mr} переменного радиуса $p$
для любого $\mu$-измеримого множества $S$ согласно \eqref{DzEk} при $d\leq 2$ имеем 
\begin{equation*}
{\mathfrak m}_d^p(E\cap S)\overset{\eqref{mr}}{\leq}  \dfrac{\pi^{d/2}}{{\Gamma(1+d}/2)}\sum_{S\cap \overline D_{z_k}(t_k)\neq \varnothing} t_{k}^d\leq \pi\sum_{S\cap \overline D_{z_k}(t_k)\neq \varnothing} t_{k}^d,
\end{equation*}
откуда по  неравенству  \eqref{SDyU}  получаем  \eqref{mmES}, и теорема \ref{lemmuq} доказана.
\end{proof}

\begin{theorem}\label{lem8_2}
Пусть функция  $r\colon \CC\to (0,1]$ удовлетворяет условию \eqref{qr},  
или эквивалентному ему условию \eqref{nubstrr} из замечания\/ {\rm \ref{rem3},}
а $\mu$ --- распределение масс конечного порядка на $\CC$.  Тогда для  любого   $d\in (0,2]$ можно подобрать функцию $p\colon \CC\to\RR^+\setminus 0$, также удовлетворяющую условию \eqref{qr}, для которой 
\begin{equation}\label{pleqr}
\inf_{z\in \CC} \frac{\ln p(z)}{\ln(2+ |z|)}>-\infty, \quad p(z)\underset{z\in \CC}{\leq} r(z), 
\end{equation}
а  множество $E$, определённое в \eqref{Er} с этой функцией $p$,
для каждого борелевского подмножества $S\subset \CC$ удовлетворяет  ограничению  
\begin{equation}\label{mmES1}
{\mathfrak m}_d^r(E\cap S) \leq \sup_{z\in S} r(z)\leq 1.
\end{equation}
\end{theorem}

\begin{proof}
Для функции $r$ с  условием \eqref{qr} в форме  \eqref{nubstrr} из замечания \ref{rem3} 
найдётся достаточно большое $Q\in \RR^+$, для которого 
\begin{equation}\label{Qd}
\frac{1}{(2+|z|)^Q}\underset{z\in \CC}{\leq}r(z).
\end{equation}
Если, как в условиях заключительной части основной теоремы, 
$$ 
\ord[\mu]<l<+\infty,
$$ 
то для $\mu$ найдётся   число $C\geq 1$, для  которого  
\begin{equation}\label{vadu}
 \mu^{\rad}(t)\leq C(1+t)^l \quad\text{при всех $t\in \RR^+$}.   
\end{equation}
Для заданного  $d\in (0,2]$ положим 
\begin{equation}\label{PQd}
P:=Q+1+\frac{1}{d}(Q+l+1)
\end{equation}
и рассмотрим отличную от \eqref{pP}  функцию  
\begin{equation}\label{pPP}
p\colon z\underset{z\in \CC}{\longmapsto} \frac{1}{(60(l+1)C)^{1/d}(4+|z|)^P}, 
\end{equation}
для которой  ввиду $C\geq 1$, $d>0$ и $P\overset{\eqref{PQd}}{\geq} 1$ имеем
\begin{equation}\label{pPPs}
s\overset{\eqref{r0}}{:=}\sup_{z\in \CC} p(z)\leq \frac{1}{2}.  
\end{equation}
Очевидно, любые функции \eqref{pPP} удовлетворяют  условию  \eqref{qr}, так как
\begin{equation}\label{pqp}
\inf_{z\in \CC}\frac{\ln p(z)}{\ln(2+|z|)}=-P>-\infty, 
\end{equation}
т.е. первому соотношению в \eqref{pleqr}. Поэтому далее можем неограниченно увеличивать $P\geq 1$, оставаясь в рамках условия \eqref{pqp} для функции $p$ из \eqref{pPP}. 

В обозначении \eqref{Sdrcup} согласно  \eqref{pPPs} имеем включение 
\begin{equation}\label{S4}
S^{\cup(2s)}\overset{\eqref{pPPs}}{\subset} S^{\cup 1}\quad\text{для любого  $S\subset \CC$,}
\end{equation}
а также  из \eqref{Qd} согласно выбору $P$ в \eqref{PQd}  ввиду  $C\geq 1$ получаем 
\begin{equation}\label{pqPe}
p(z)\overset{\eqref{pPP},\eqref{PQd}}{\underset{z\in \CC}{\leq}}
\frac{1}{(2+|z|)^Q} \underset{z\in \CC}{\overset{\eqref{Qd}}{\leq}} r(z),
\end{equation}
что даёт всё \eqref{pleqr}. По теореме \ref{lemmuq}  для множества $E$, определённого  в \eqref{Er}, 
для каждого подмножества $S\subset \CC$ выполнено неравенство \eqref{mmES}, где 
\begin{equation}\label{Sss}
S_s\subset S^{\cup(2s)}\overset{\eqref{S4}}{\subset} S^{\cup 1}.  
\end{equation}
Кроме того, в обозначениях \eqref{Sdrp} ввиду  $d\in (0,2]$ получаем 
\begin{equation}\label{pdP}
(p^d)^{\vee s}(z) \overset{\eqref{pPPs}}{\leq} (p^d)^{\vee \frac12}(z) 
\overset{\eqref{pPP},\eqref{PQd}}{\underset{z\in \CC}{\leq}} 
\frac{1}{60C(l+1)(3+|z|)^{Q+l+1}}.
\end{equation}
откуда согласно \eqref{vadu}  для интеграла с множителем $60$ в правой части \eqref{mmES} с  распределением  масс 
 $\mu$  имеем 
\begin{multline}\label{mmE}
 60 \int_{S_s} \bigl(p^d\bigr)^{\vee s} \dd \mu
\overset{\eqref{pdP}}{\leq} 60
\int_{S_s} \frac{\dd \mu(z)}{60C(l+1)(3+|z|)^{Q+l+1}}
\\
\overset{\eqref{Sss}}{\leq} \sup_{z\in S^{\cup 1}} \frac{1}{(3+|z|)^Q}
\int_{\CC} \frac{\dd \mu(z)}{C(l+1)(3+|z|)^{l+1}}
\\
=\frac{1}{\Bigl(3+\inf\limits_{z\in S} |z|-1\Bigr)^Q}
\int_0^{+\infty} \frac{\dd \mu^{\rad}(t)}{C(l+1)(3+t)^{l+1}}\\
\overset{\eqref{vadu}}{\leq} 
\sup_{z\in S} \frac{1}{(2+|z|)^Q}\int_0^{+\infty}\frac{C(1+t)^l\dd t}{C(3+t)^{l+2}}
\leq \sup_{z\in S} \frac{1}{(2+|z|)^Q}\leq \sup_{z\in S} r(z).
\end{multline}
Таким образом, имеют место неравенства  
\begin{equation}\label{mmESr}
{\mathfrak m}_d^p(E\cap S) \overset{\eqref{mmES}}{\leq}  60 \int_{S_s} \bigl(p^d\bigr)^{\vee s} \dd \mu \overset{\eqref{mmE}}{\leq} 
 \sup_{z\in S} r(z).
\end{equation}
Но по построению $r\geq p$ всюду, вследствие чего 
${\mathfrak m}_d^r(E\cap S)\overset{\eqref{mts}}{\leq} {\mathfrak m}_d^p(E\cap S)$.  
Отсюда согласно  \eqref{mmESr} получаем \eqref{mmES1}, и теорема \ref{lem8_2} доказана. 
\end{proof}

Вернёмся к доказательству заключительной части основной	теоремы  для случай распределения масс Рисса $\mu:=\varDelta_u$
конечного порядка 
$$
\ord[\mu]\overset{\eqref{ordnu}}{=}\ord[\varDelta_u]<+\infty.
$$ 

По теореме \ref{lem8_2} можно построить   функцию $p\colon \CC\to (0,1]$    со свойствами \eqref{pleqr} и  \eqref{mmES1}, удовлетворяющую условию \eqref{qr}, т.е. первому соотношению в \eqref{pleqr}. Отсюда по уже доказанной части основной теоремы  для этой функции $p$ в роли $r$ существует целая функция $f_p\not\equiv 0$, удовлетворяющая неравенствам 
\begin{equation}\label{ubfp}
\ln |f_p(z)|\underset{z\in \CC}{\overset{\eqref{ubc}}{\leq}} u^{\bullet p}(z) \underset{z\in \CC}{\overset{\eqref{ubc}}{\leq}} u^{\circ p}(z).
\end{equation}
 По  формуле Пуассона\,--\,Йенсена\,--\,При\-в\-а\-л\-о\-ва 
 \cite[гл. II, \S~2]{Privalov}, \cite[4.5]{Rans}, \cite[3.7]{HK}
для распределения масс $\mu:=\varDelta_u$ имеем 
\begin{equation*}
u^{\circ p}(z)= u(z)+\int_0^{p(z)}\frac{\mu^{\rad}(t)}{t}\dd t
\quad\text{при всех  $z\overset{\eqref{fytu}}{\in} \mathbb C\setminus (-\infty)_u$}. 
\end{equation*}
Отсюда по определению \eqref{Er} множества $E$ из теоремы \ref{lem8_2}  по \eqref{ubfp} получаем
\begin{equation}\label{ubfpd}
\ln |f_p(z)|\overset{\eqref{ubfp}}{\leq} u(z)+\frac1d \quad\text{при всех $z\in \CC\setminus \bigl(E\cup (-\infty)_u\bigr)$,}
\end{equation}
где $(-\infty)_u\subset E$ и $E\cup (-\infty)_u=E$, поскольку 
$$
\int_0^{p(z)}\frac{\mu^{\rad}(t)}{t}\dd t=+\infty \quad\text{для всех  $z\in (-\infty)_u$}. 
$$
Положим $f:=f_pe^{-1/d}\not\equiv 0$. Тогда из  \eqref{ubfpd} имеем  \eqref{{uE1}u}, а именно:
\begin{equation}\label{ubfpd+}
\ln |f(z)|\underset{z\in \CC}{\equiv}\ln |f_p(z)|-\frac{1}{d}\overset{\eqref{ubfpd}}{\leq} u(z)\quad\text{при всех $z\in \CC\setminus E$,}
\end{equation}
где по теореме \ref{lem8_2} для множества $E$ выполнено \eqref{mmES1}, или \eqref{{uE1}E}, для любого подмножества  $S\subset \CC$. 
В то же время для функции  $\ln|f|=\ln|f_p|-1/d$  по \eqref{ubfp} и неравенству 
$p\overset{\eqref{pleqr}}{\leq} r$ на $\CC$ по-прежнему выполнены и соотношения из \eqref{ubc}:
\begin{equation*}\label{ubfpp}
\ln |f(z)|\overset{\eqref{ubfpd+}}{\underset{z\in \CC}{\leq}}\ln |f_p(z)|\underset{z\in \CC}{\overset{\eqref{ubfp}}{\leq}} u^{\bullet p}(z) \underset{z\in \CC}{\overset{\eqref{pleqr}}{\leq}} u^{\bullet r}(z)
\underset{z\in \CC}{\overset{\eqref{vbc}}{\leq}} u^{\circ r}(z)
\underset{z\in \CC}{\leq} \sup_{D_z(r(z))}u,
\end{equation*}
и, как следствие, \eqref{Mulnf}.  Это завершает доказательство основной теоремы.

\end{fulltext}

\end{document}